\documentclass{sig-alternate-per}
\usepackage{theorem}
\usepackage{ifpdf}
\usepackage{hyperref}
\usepackage{setspace}

\hypersetup{
    bookmarksopen=false,
    bookmarksnumbered=true,
}
\ifpdf
  \hypersetup{colorlinks=true,linkcolor=black,urlcolor=black,citecolor=black,pdfpagemode=UseOutlines,plainpages=false,pdfpagelabels}
\else
  \hypersetup{colorlinks=true,linkcolor=black,urlcolor=black,citecolor=black,pdfpagemode=UseOutlines,plainpages=false,pdfpagelabels}
\fi

\newtheorem{theorem}{Theorem}

\newtheorem{remark}{Remark}
\newtheorem{lemma}{Lemma}

\newcommand{\ii}{i=1,\dots,N}
\newcommand{\jj}{j=1,\dots,N}

\newcommand{\ul}{\underline}
\newcommand{\beq}{\begin{equation}}
\newcommand{\eeq}{\end{equation}}
\newcommand{\btilde}{\tilde{b}}
\newcommand{\Btilde}{\widetilde{B}}
\newcommand{\E}{\mathbb{E}}
\newcommand{\limdist}{\,{\buildrel d \over \rightarrow}\,}

\begin{document}

\conferenceinfo{Performance}{2011 Amsterdam, The Netherlands}

\title{Queueing networks with a single shared server:\\Light and heavy traffic}

\numberofauthors{3}
\author{
\alignauthor
M.A.A. Boon \\
       \affaddr{Department of Mathematics and Computer Science}\\
       \affaddr{Eindhoven University of Technology}\\
       \affaddr{ P.O. Box 513, 5600MB Eindhoven, The Netherlands}\\
       \email{marko@win.tue.nl}
\alignauthor
R.D. van der Mei \\
       \affaddr{Centre for Mathematics and Computer Science (CWI)}\\
       \affaddr{Department of Probability and Stochastic Networks}\\
       \affaddr{1098 SJ Amsterdam, The Netherlands}\\
       \email{mei@cwi.nl}
\alignauthor
E.M.M. Winands \\
       \affaddr{Department of Mathematics, Section Stochastics}\\
       \affaddr{VU University}\\
       \affaddr{De Boelelaan 1081a, 1081HV Amsterdam, The Netherlands}\\
       \email{emm.winands@few.vu.nl}
}

\date{30 July 1999}

\maketitle
\begin{abstract}
We study a queueing network with a single shared server, that serves the queues in a cyclic order according to the gated service discipline.
External customers arrive at the queues according to independent Poisson processes. After completing service, a customer either leaves the system or is routed to another queue. This model is very generic and finds many applications in computer systems, communication networks, manufacturing
systems and robotics. Special cases of the introduced network include well-known polling models and tandem queues. We derive exact limits of the mean delays under both heavy-traffic and light-traffic conditions. By interpolating between these asymptotic regimes,  we develop simple closed-form approximations for the mean delays for arbitrary loads.
\end{abstract}




\section{Introduction}

In this paper we study a queueing network served by a single shared server, that visits the queues in a cyclic order. Each queue receives gated service. Customers from the outside arrive at the queues according to independent Poisson processes, and the service time and switch-over time
distributions are general. After receiving service at queue $i$, a customer is either routed to queue $j$ with probability $p_{i,j}$, or leaves the system with probability $p_{i,0}$. This model can be seen as an extension of the standard polling model by customer routing, introduced and analyzed by Sidi et al. \cite{sidi1,sidi2}.
The possibility of re-routing of customers further enhances the already-extensive modeling capabilities of polling models,
which find applications in diverse areas such as computer systems, communication networks, logistics, flexible manufacturing
systems, robotics systems, production systems and maintenance systems (see, for example, \cite{boonapplications2011} for overviews). Applications of this type of customer routing can be found, for example, in manufacturing when products undergo service in a number of stages or in the context of rework. The present research can be seen as a unifying analysis for a variety of special cases (besides polling systems), such as tandem queues \cite{nair}, multi-stage queueing models with parallel servers \cite{katayama}, feedback vacation queues \cite{boxmayechiali97} and many others.

The main contribution of the present paper is twofold. Firstly,  we study exact  heavy-traffic (HT) asymptotics of the system under consideration. Under HT scalings we derive a closed-form expression for the joint queue-length vector at polling instants under HT scalings. This result, in turn, is shown to lead to a closed-form expression for the expected delay at a queue at an arbitrary moment. This expression is strikingly simple and shows explicitly how the  expected delays depend on the system parameters, and in particular, on the routing probabilities $p_{i,j}$. Secondly, we derive  a closed-form approximation for the mean delay for arbitrary loads, based on  an interpolation between the light-traffic (LT) and HT limits. Numerical results are presented to assess the accuracy of this approximation. We would like to note that the analysis of the present paper can be extended in several directions, for example different server routing policies or service disciplines. For reasons of compactness, these results are however not discussed.

The remainder of this paper is organized as follows. In Section \ref{modelsection} we describe the model in further detail. This model is analyzed under HT scalings in Section \ref{analysissection}. In Section \ref{approxsection} we develop an approximation based on the LT and HT limits for the mean waiting times and give a numerical example.

\section{Model description}\label{modelsection}

In this paper we consider a queueing network consisting of $N\geq2$ infinite buffer queues $Q_1,\dots,Q_N$.
External customers arrive at $Q_i$ according to a Poisson arrival process with rate $\lambda_i$, and have a generally distributed service requirement $B_i$ at $Q_i$, with mean value $b_i := E[B_i]$, and second moment $b_i^{(2)} := E[B_i^2]$. The queues are served by a single server in cyclic order. Whenever the server switches from $Q_i$ to $Q_{i+1}$, a switch-over time $R_i$ is incurred, with mean $r_i$. The cycle time $C_i$ is the time between successive moments when the server arrives at $Q_i$. The total switch-over time in a cycle is denoted by $R=\sum_{i=1}^N R_i$ 
and its first two moments are $r:=\E[R]$ and $r^{(2)}:=\E[R^2]$. Indices throughout the paper are modulo $N$, so $Q_{N+1}$ actually refers to $Q_1$.
All service times and switch-over times are mutually independent. Each queue receives gated service, which means that only those customers present at the server's arrival at $Q_i$ will be served before the server switches to the next queue. This queueing network can be modeled as a \emph{polling system} with the specific feature that it allows for routing of the customers: upon completion of service at $Q_i$, a customer is either routed to $Q_j$ with probability $p_{i,j}$, or leaves the system with probability $p_{i,0}$. Note that $\sum_{j=0}^N~p_{i,j}=1$ for all $i$, and that the transition of a customer from $Q_i$ to $Q_j$ takes no time. The model under consideration has a branching structure, which is discussed in more detail by Resing \cite{resing93}.  The total arrival rate at $Q_i$ is denoted by $\gamma_i$, which is the unique solution of the following set of linear equations:
\begin{equation*}
\gamma_i = \lambda_i + \sum_{j=1}^N \gamma_j p_{j,i},\qquad\ii.
\end{equation*}
The offered load to $Q_i$ is $\rho_i:=\gamma_i b_i$ and the total utilisation is $\rho:=\sum_{i=1}^N \rho_i$. We assume that the system is stable, which means that $\rho$ should be less than one (see \cite{sidi2}). The total service time of a customer is the total amount of service given during the presence
of the customer in the network, denoted by $\Btilde_i$, and its first two moments by
$\btilde_i := E[\Btilde_i]$ and $\btilde_i^{(2)} := E[\Btilde_i^2]$.
The first two moments are uniquely determined by the following set of linear equations: For $\ii$,
\begin{eqnarray*}
\btilde_i &=& b_i + \sum_{j=1}^N \btilde_j p_{i,j},\\
\btilde_i^{(2)} &=& b_i^{(2)} + 2b_i\sum_{j=1}^N \btilde_j p_{i,j}+ \sum_{j=1}^N \btilde_j^{(2)} p_{i,j}.
\end{eqnarray*}

For each variable $x$ (which may be a scalar, a vector or a matrix) that is a function of $\rho$, we denote its value {\it evaluated at} $\rho=1$ by $\hat{x}$. Also, we will be taking HT limits, letting $\rho \uparrow 1$. To be precise, the limit is taken such that the arrival rates $\lambda_1,\ldots,\lambda_N$ are increased, while keeping the service and switch-over time distributions, the routing
probabilities and the {\it ratios} between these arrival rates fixed. Also, an $N$-dimensional vector ${\ul x}$ has components $(x_1,\ldots,x_N)$.

\section{Heavy-traffic analysis}\label{analysissection}

Let ${\ul X}:=\left(X_1 \ldots, X_N \right)$ be the $N$-dimensional vector that describes the joint queue length at a visit beginning at $Q_1$.
Sidi et al. \cite{sidi2} show that the joint queue-length process at successive visit beginnings to $Q_1$ constitutes an $N$-dimensional multi-type branching process (MTBP) with immigration in each state.
In this section we focus on the limiting behavior of ${\ul X}$ as $\rho$ goes to 1.
\noindent\begin{theorem}
The joint queue-length vector at polling instants at $Q_1$ has the following asymptotic behavior:
\[
(1-\rho)
\left(
\begin{array}{c}
X_1 \\ \vdots \\ X_N
\end{array}
\right)
\limdist
{\tilde{b}^{(2)} \over 2\tilde{b}^{(1)}}
~
{1 \over \delta}
\left(
\begin{array}{c}
\hat{u}_1 \\ \vdots \\ \hat{u}_N
\end{array}
\right)
~
\Gamma(\alpha, 1)~~~(\rho \uparrow 1),
\]
with
\begin{align}
{u}_i &:={\lambda}_i \sum_{j=i}^N {\rho}_j + \sum_{j=i}^N \gamma_j p_{j,i}~~~(\ii),\nonumber\\
\delta &:=  \hat{\ul u}^{\top} \ul{\tilde{b}}
=
\sum_{i=1}^N   \hat{\lambda}_i \tilde{b}_i \sum_{j=i+1}^N  \hat{\rho}_j+ \sum_{i=1}^N \tilde{b}_i \sum_{j=i}^N \hat{\gamma}_j p_{j,i},\label{delta}\\
\alpha &:= 2 r \delta {\tilde{b}^{(1)} \over \tilde{b}^{(2)}},\nonumber\\
\tilde{b}^{(k)}&:=\sum_{i=1}^N\lambda_i\E[\Btilde_i^k]/\sum_{j=1}^N\lambda_j,\nonumber
\end{align}
and where $\Gamma(\alpha, 1)$ is a gamma-distributed random variable with shape parameter $\alpha$ and scale parameter 1.
\end{theorem}
For reasons of compactness we omit the proof, noting that the basis for the
proof is given by the general framework as proposed in \cite{RvdM_QUESTA}.

Next, we focus on the workload in the individual queues. To this end, we note that simple balancing arguments can be used to show that $\E[C_i]=r/(1-\rho)$, which does not depend on $i$.
To obtain HT-results for the amount of work in each queue, we use the Heavy Traffic Averaging Principle (HTAP) for polling systems \cite{coffman98}. When the system becomes saturated, two limiting processes take place. The scaled total workload tends to a Bessel-type diffusion whereas the work in each queue is changing at a much faster rate than the total workload. This implies that during the course of a cycle, the total workload can be considered as constant, while the workloads of the individual queues fluctuate according to a fluid model. The HTAP relates these two limiting processes. The fluid limit of the per-queue workload is obtained by dividing by $r/(1-\rho)$ and letting $\rho \uparrow 1$. For our model, the fluid model for the workload at $Q_i$ is a piecewise linear function. More precisely, it is easy to show that  the fluid limit of the mean amount of work at $Q_i$ at the beginning of  a visit to $Q_j$ is $\sum_{k=i}^{j+N-1} \hat{\gamma}_k \tilde{b}_i(\hat{\lambda}_i b_k + p_{k,i})$ for $j \neq i$ and $\hat{\gamma}_i b_i$ for $j=i$. Moreover, in the fluid limit the probability that at an arbitrary moment the server is visiting $Q_j$ is $\hat{\rho}_j~(\jj)$. Combining these observations, one can obtain the following expression for $\delta_i$, defined as the fluid limit of the average amount of work at $Q_i$.

\hspace*{-2cm}\noindent\begin{lemma}
For $\ii$,
\begin{align}
\delta_i&
=
{1 \over 2} \hat{\rho}_i \hat{\gamma}_i \tilde{b}_i ( 1 + \hat{\lambda}_i b_i + p_{i,i})
\\
&+
\sum_{j=i+1}^{i+N-1} \hat{\rho}_j
\big(
{1 \over 2} \hat{\gamma}_j \tilde{b}_i(\hat{\lambda}_i \tilde{b}_j + p_{j,i})
+
\sum_{k=i}^{j-1} \hat{\gamma}_k \tilde{b}_i ( \hat{\lambda}_i b_k + p_{k,i})
\big).\nonumber
\end{align}
\end{lemma}
Note that in the classical case where $p_{i,j}=0$ for all $i,j$ we have $\delta_i = \hat{\rho}_i(1+\hat{\rho}_i)/2$. Moreover, it is easily verified that $\sum_{i=1}^N \delta_i=\delta$, where $\delta$ is defined in \eqref{delta}.

Subsequently, the diffusion limit of the \emph{total} workload process and the workload in the individual queues can be related using the HTAP. To this end, we start with the cycle-time distribution under HT scalings.

\noindent\begin{lemma}\label{lemmacycletime}
For $\ii$,
\[
(1-\rho)C_i \limdist \Gamma(\alpha,\mu),\qquad(\rho\uparrow1),
\]
with $\alpha=2r\delta \btilde^{(1)}/\btilde^{(2)}$ and $\mu=2\delta\btilde^{(1)}/\btilde^{(2)}$.
\end{lemma}
Note that the distribution of $C_i$ no longer depends on $i$ in the HT limit. The proof can be found along the same lines as in \cite{RvdM_QUESTA}. Lemma~\ref{lemmacycletime} implies that the HT limit of the \emph{mean} (scaled) amount of work found by an arbitrary customer is $\delta_i\E[\mathbf{C}]$, where $\mathbf{C}$ is the \emph{length-biased} cycle time, with limiting distribution $\Gamma(\alpha+1,\mu)$. Again, see \cite{RvdM_QUESTA} for more details.

The (HT limit of the) mean waiting time of an arbitrary customer in $Q_i$ can be found by application of Little's Law to the mean queue length at $Q_i$, which is simply the mean amount of work in $Q_i$ divided by the mean total service time.
\noindent\begin{theorem}
For $\ii$,
\begin{equation}
(1-\rho)\E[W_i] \rightarrow \left(r+\frac{\btilde^{(2)}}{2\delta\btilde^{(1)}}\right)\frac{\delta_i}{\btilde_i\hat\gamma_i}, \qquad (\rho\uparrow1).
\label{EWht}
\end{equation}
\end{theorem}

\begin{remark}[Insensitivity]
Equation \eqref{EWht} reveals a variety of properties about the dependence of the limiting mean delay with respect to the system parameters.
The mean waiting times $\E[W_i]$ are independent of the visit order of the server, depend on the switch-over time distributions only through $r$, and depend on the service-time distributions only through $\btilde^{(1)}$ and $\btilde^{(2)}$.
\end{remark}

\section{Approximation}\label{approxsection}

The LT limit of $\E[W_i]$ can be found by conditioning on the customer type (external or internally routed).
\noindent\begin{theorem}
For $\ii$,
\begin{equation}
\E[W_i] \rightarrow \frac{\lambda_i}{\gamma_i}\frac{r^{(2)}}{2r} + \sum_{j=i-N}^{i-1}\frac{\gamma_jp_{j,i}}{\gamma_i}\sum_{k=j}^{i-1}r_k, \qquad (\rho\downarrow0).
\label{EWlt}
\end{equation}
\end{theorem}
In light traffic we ignore all $O(\rho)$ terms, which implies that we can consider a customer as being alone in the system. Equation \eqref{EWlt} can be interpreted as follows. An arbitrary customer in $Q_i$ has arrived from outside the network with probability $\frac{\lambda_i}{\gamma_i}$. In this case he has to wait for a residual total switchover time with mean ${r^{(2)}}/{2r}$. If a customer in $Q_i$ arrives after being served in another queue, say $Q_j$ (with probability ${\gamma_jp_{j,i}}/{\gamma_i}$), he has to wait for the mean switch-over times $r_j,\dots,r_{i-1}$.

Subsequently, we construct an interpolation between the LT and HT limits that can be used as an approximation for the mean waiting times. For $\ii$,
\begin{equation}
W_i^\textit{approx}=\frac{w_i^\textit{LT} + (w_i^\textit{HT}-w_i^\textit{LT})\rho}{1-\rho},\label{EWapprox}
\end{equation}
where $w_i^\textit{LT}$ and $w_i^\textit{HT}$ are the LT and HT limits respectively, as given in \eqref{EWlt} and \eqref{EWht}. Because of the way $W_i^\textit{approx}$ is constructed, it has the nice properties that it is exact as $\rho\downarrow0$ and $\rho\uparrow1$. Furthermore, it satisfies a so-called pseudo-conservation law for the mean waiting times, which is derived in \cite{sidi2}. This implies that the $W_i^\textit{approx}$ yields exact results for symmetric (and, hence, single-queue) systems. 

We do not aim at giving an extensive numerical study to assess the accuracy of the approximation. Instead, we give one numerical example that indicates the versatility of the model that we have discussed, and shows the practical usage of the approximation \eqref{EWapprox}. To this end, we use an example that was introduced by Katayama \cite{katayama}, who studies a network consisting of three queues. Customers arrive at $Q_1$ and $Q_2$, and are routed to $Q_3$ after being served. This model, which is referred to as a tandem queueing model with parallel queues in the first stage, is a special case of the model discussed in the present paper. We simply put $p_{1,3}=p_{2,3}=p_{3,0}=1$ and all other $p_{i,j}$ are zero. We use the same values as in \cite{katayama}: $\lambda_1=\lambda_2/10$, service times are deterministic with $b_1=b_2=1$, and $b_3=5$. The server visits the queues in cyclic order: 1, 2, 3, 1, \dots. The only difference with the model discussed in \cite{katayama} is that we introduce (deterministic) switch-over times $r_2=r_3=2$. We assume that no time is required to switch between the two queues in the first stage, so $r_1=0$. In Table \ref{numericalresults} we show the mean waiting times of customers at the three queues and their approximated values. From this table we can see that the accuracy is best for values of $\rho$ close to 0 or 1, but the overall accuracy is very good in general.

\begin{table}[h]
\[
\begin{array}{|l|ccccccc|}
\hline
\rho & 0.01 & 0.1 & 0.3 & 0.5 & 0.7 & 0.9 & 0.99 \\
\hline
\E[W_1]             &  2.05 & 2.53 & 3.87 & 6.07 & 10.95 & 34.87 & 356.60 \\
W_1^\textit{approx} &  2.04 & 2.40 & 3.53 & 5.57 & 10.34 & 34.17 & 355.86 \\
\hline
\hline
\E[W_2] &             2.05 & 2.56 & 4.01 & 6.45 & 11.95 & 39.05 & 403.97 \\
W_2^\textit{approx}&   2.04 & 2.45 & 3.74 & 6.05 & 11.46 & 38.49 & 403.39\\
\hline
\hline
\E[W_3] &              2.02 & 2.26 & 3.18 & 5.04 & 9.62 & 33.00 & 349.85 \\
W_3^\textit{approx} &  2.04 & 2.39 & 3.51 & 5.52 & 10.22 & 33.69 & 350.57\\
\hline
\end{array}
\]
\caption{Results for the numerical example.}
\label{numericalresults}
\end{table}

\bibliographystyle{abbrv}

\end{document}